\title{Hamilton cycles in highly connected and expanding graphs}
\author{{Dan Hefetz \thanks{School of Computer Science,
Raymond and Beverly Sackler Faculty of Exact Sciences, Tel Aviv
University, Tel Aviv, 69978, Israel. Email: dannyh@post.tau.ac.il.
This paper is a part of the author's Ph.D. thesis written under the
supervision of Prof. Michael Krivelevich.}} \quad {Michael
Krivelevich \thanks{ Department of Mathematics, Raymond and Beverly
Sackler Faculty of Exact Sciences, Tel Aviv University, Tel Aviv,
69978, Israel. Email: krivelev@post.tau.ac.il. Research supported in
part by a USA-Israeli BSF grant and a grant from the Israeli Science
Foundation.}} \quad {Tibor Szab\'o \thanks{Institute of Theoretical
Computer Science, ETH Zurich, CH-8092 Switzerland. Email:
szabo@inf.ethz.ch.}}}

\documentclass[12pt]{article}

\usepackage{amsmath,amssymb,latexsym,color,epsfig}

\newif\ifnotesw\noteswtrue

\def\a{\alpha}

\def\r{\rho}

\def\s{\sigma}

\newtheorem{theorem}{Theorem}[section]

\newtheorem{lemma}[theorem]{Lemma}

\newtheorem{claim}[theorem]{Claim}

\newtheorem{prop}[theorem]{Proposition}

\newcommand{\proofstart}{{\bf Proof\hspace{2em}}}

\newcommand{\proofend}{\hspace*{\fill}\mbox{$\Box$}}

\newcommand{\rdown}[1]{{\mbox{$ \lfloor #1 \rfloor $}}}

\renewcommand{\epsilon}{\varepsilon}

\begin{document}

\maketitle

\begin{abstract}
In this paper we prove a sufficient condition for the existence of a
Hamilton cycle, which is applicable to a wide variety of graphs,
including relatively sparse graphs. In contrast to previous
criteria, ours is based on only two properties: one requiring
expansion of ``small'' sets, the other ensuring the existence of an
edge between any two disjoint ``large'' sets. We also discuss
applications in positional games, random graphs and extremal graph
theory.
\end{abstract}

\section{Introduction} \label{sec::intro}
A Hamilton cycle in a graph $G$ is a cycle passing through all
vertices of $G$.  A graph is called \emph{Hamiltonian} if it admits
a Hamilton cycle. Hamiltonicity is one of the most central notions
in Graph Theory, and many efforts have been devoted to obtain
sufficient conditions for the existence of a Hamilton cycle (a
"nice" necessary and sufficient condition should not be expected
however, as deciding whether a given graph contains a Hamilton cycle
is known to be NP-complete). In this paper we will mostly concern
ourselves with establishing a sufficient condition for Hamiltonicity
which is applicable to a wide class of sparse graphs.

One of the first Hamiltonicity results is
the celebrated theorem of Dirac~\cite{Dirac}, which asserts that if
the minimum degree of a graph $G$ on $n$ vertices is at least $n/2$
then $G$ is Hamiltonian. Since then, many other sufficient
conditions that deal with dense graphs, were obtained (see
e.g.~\cite{Gould} for a comprehensive reference). However, all these
conditions require the graph to have $\Theta(n^2)$ edges whereas for
a Hamilton cycle, only $n$ edges are needed. Chv\'atal and
Erd\H{o}s~\cite{CE} proved that if $\kappa(G) \geq \alpha(G)$ (that
is, the vertex connectivity of $G$ is at least as large as the size
of a largest independent set in $G$) then $G$ is Hamiltonian. Note
that if $G$ is a $d$-regular graph, then $\kappa(G) \leq d$ and
$\alpha(G) \geq \frac{n}{d+1}$; hence the Chv\'atal-Erd\H{o}s
criterion cannot be applied if $d \leq c \sqrt{n}$ for an
appropriate constant $c$.







When looking for sufficient conditions for the Hamiltonicity of
sparse graphs, it is natural to look at random graphs with an
appropriate edge probability. Erd\H{o}s and R\'enyi~\cite{ER61}
raised the question of what is the threshold for Hamiltonicity in
random graphs. After a series of efforts by various researchers,
including Korshunov~\cite{Korsh76} and P\'osa~\cite{Posa}, the
problem was finally solved by Koml\'os and Szemer\'edi~\cite{KS},
who proved that if $p=(\log n + \log \log n +\omega(1))/n$, where
$\omega(1)$ tends to infinity with $n$ arbitrarily slowly, then
$G(n,p)$ is a.s. Hamiltonian. Note that this is best possible since
for $p \leq (\log n + \log \log n - \omega(1))/n$ almost surely
there are vertices of degree at most one in $G(n,p)$.

The next natural step is to look for Hamilton cycles in relatively
sparse pseudo-random graphs. During the last few years, several such
sufficient conditions were found (see e.g.~\cite{FK, KrSu}). These
are quite complicated at times as they rely on many properties of
pseudo-random graphs. Furthermore, one can argue that these
conditions are not the most natural, as Hamiltonicity is a monotone
increasing property, whereas pseudo-randomness is not. Our main
result is a natural and simple (at least on the qualitative level)
sufficient condition based on expansion and high connectivity.
Before stating the result we introduce and discuss the following
properties of a graph $G=(V,E)$ where $|V|=n$. As usual, the
notation $N(S)$ stands for the {\em external neighborhood} of $S$,
that is, $N(S)=\{v \in V \setminus S : \exists u \in S, \; (u,v) \in
E\}$. Let $d=d(n)$ be a parameter.

\begin{description}
\item[P1] For every $S \subset V$, if $|S|\le \frac{n \log \log n \log d}{d\log n \log \log \log n}$
then $|N(S)|\ge d|S|$;
\item[P2] There is an edge in $G$ between any two disjoint
subsets $A,B \subseteq V$ such that $|A|,|B|\ge \frac{n \log \log n
\log d}{4130 \log n \log \log \log n}$.
\end{description}

From now on, for the sake of convenience, we denote
$$
m=m(n,d)=\frac{\log n\cdot\log\log\log n}{\log\log n\cdot\log d}\ .
$$

Let us give an informal interpretation of the above conditions.
Condition P1 guarantees {\em expansion}: every sufficiently small
vertex subset (of size $|S|\le \frac{n}{dm}$) expands by a factor of
$d$. Condition P2 is what can be classified as a {\em high
connectivity} condition of some sort: every two disjoint subsets $A,
B \subseteq V$ which are relatively large (of size $|A|, |B| \geq
\frac{n}{4130m}$) are connected by at least one edge. Note that
properties P1 and P$2$ together guarantee some expansion for {\em
every} $S \subset V(G)$ of size $o(n)$. Indeed, if $|S|\le
\frac{n}{d m}$ then $|N(S)|\ge d|S|$ by property P1. If $\frac{n}{d
m} < |S| < \frac{n}{4130 m}$ (assuming $d > 4130$) then $S$ contains
a subset of size exactly $\frac{n}{d m}$ and so by property P1
expands at least to a size of $\frac{n}{m}$, that is it expands by a
factor of at least 4130. Finally, if $|S| \geq \frac{n}{4130 m}$
then $N(S) \geq (1-o(1))n$ as, by property P$2$, the number of
vertices of $V \setminus S$ that do not have any neighbor in $S$ is
strictly less than $\frac{n}{4130 m}$.

We can now state our main result:

\begin{theorem}\label{th}
Let $12 \leq d \leq e^{\sqrt[3]{\log n}}$ and let $G$ be a graph on
$n$ vertices satisfying properties P1, P$2$ as above; then $G$ is
Hamiltonian, for sufficiently large $n$.
\end{theorem}

The lower bound on $d$ in the theorem above can probably be somewhat
improved through a more careful implementation of our arguments. As
for the upper bound $d\le e^{\sqrt[3]{\log n}}$, it is a mere
technicality, as one expects that proving that denser graphs (that
is, graphs for which $d$ is larger) are Hamiltonian should in fact
be easier. The requirement $d \leq e^{\sqrt[3]{\log n}}$ makes sure
(in particular) that $\frac{n}{4130 m} = o(n)$ and so P2 is a
non-trivial condition. We can obtain a sufficient condition for
Hamiltonicity, similar to that of Theorem~\ref{th}, and applicable
to graphs with larger values of $d=d(n)$ as well; more details are
given in Section~\ref{large-d}.





%













It is instructive to observe that neither P1 nor P2 is enough to
guarantee Hamiltonicity by itself, without relying on its companion
property (unless of course they degenerate to something trivial).
Indeed, for property P1 observe that the complete graph $K_{n,n+1}$
is a very strong expander locally, yet it obviously does not contain
a Hamilton cycle. As for property P2, the graph $G$ formed by a
disjoint union of a clique of size $n-\frac{n}{4130m}+1$ and
$\frac{n}{4130m}-1$ isolated vertices clearly meets P2, but is
obviously quite far from being Hamiltonian. Thus, P1 and P2
complement each other in an essential way.




Next, we discuss several applications of our main result. Theorem
\ref{th} was first used by the authors (see~\cite{HKS}) to address  a
problem of Beck \cite{Beckbook}: they proved that
Enforcer can win the $(1,q)$ Avoider-Enforcer Hamilton cycle game,
played on the edges of $K_n$, for every $q \leq \frac{c n \log \log
\log \log n}{\log n \log \log \log n}$ where $c$ is an appropriate
constant; this is presently the best known bound. A similar result
can be obtained for Maker in the corresponding Maker-Breaker game.
The latter result falls short of the true value of the critical bias
of $\frac{c n}{\log n}$ obtained by Beck (see~\cite{Beck}), but our
proof is conceptually simpler and shorter. [A brief background: both
Maker-Breaker and Avoider-Enforcer games mentioned above are played
on the edge set of the complete graph $K_n$. In every move, Maker
(resp. Avoider) claims one unoccupied edge, then Breaker (resp.
Enforcer) responds by claiming $q$ unoccupied edges. The game ends
when all edges have been claimed by one of the players. In the
Maker-Breaker Hamiltonicity game Maker wins if he creates a Hamilton
cycle, otherwise Breaker wins. In the Avoider-Enforcer version,
Avoider wins if he avoids creating a Hamilton cycle by the end of
the game, otherwise Enforcer wins.] More details can be found
in~\cite{HKS}.

In this paper we prove several other corollaries of Theorem~\ref{th}.

A graph $G=(V,E)$ is called \emph{Hamilton-connected} if for every
$u,v \in V$ there is a Hamilton path in $G$ from $u$ to $v$.

\begin{theorem} \label{uvpath}  Let $G=(V,E)$ be a graph that satisfies properties
P1 and P2; then $G$ is Hamilton-connected.
\end{theorem}

\textbf{Remark.} An immediate consequence of Theorem~\ref{uvpath} is
that for every edge $e \in E$ there is a Hamilton cycle of $G$ that
includes $e$.

A graph $G$ is called \emph{pancyclic} if it admits a cycle of
length $k$ for every $3 \leq k \leq n$. We prove that a graph which
satisfies property P2 is "almost pancyclic".

\begin{theorem} \label{pancyclicP2}
Let $G=(V,E)$, where $|V|=n$ is sufficiently large, be a graph,
satisfying property P$2$; more precisely, for every disjoint subsets
$A,B \subseteq V$ such that $|A|,|B| \geq n/t$, where $t=t(n) \geq
2$, there is an edge between a vertex of $A$ and a vertex of $B$.
Then $G$ admits a cycle of length exactly $k$ for every $\frac{8 n
\log n}{t \log \log n} \leq k \leq n - 3n/t$.
\end{theorem}

\textbf{Remark.} The upper bound on $k$ in Theorem~\ref{pancyclicP2}
is tight up to a constant factor  in the second order term, as shown
by a disjoint union of $K_{n+1-n/t}$ and $n/t - 1$ isolated
vertices. On the other hand, we believe that the lower bound can be
improved to $\frac{c \log n}{\log t}$ for some constant $c$. Methods
recently utilized by Verstra\"ete~\cite{V} and by Sudakov and
Verstra\"ete~\cite{SV} can possibly be used to establish this
conjecture.

Theorem~\ref{th} (with minor changes to the proof) can be used to
prove the following classic result (see~\cite{KS}).

\begin{theorem} \label{randomgraph}
$G(n,p)$, where $p = (\log n + \log \log n + \omega(1))/n$, is a.s.
Hamiltonian.
\end{theorem}

Let $G=(V,E)$, where $|V|=n$, and let $f : \mathbb{Z}^+ \rightarrow
\mathbb{R}$. A pair $(A,B)$ of proper subsets of $V$ is called a
\emph{separation} of $G$ if $A \cup B = V$ and there are no edges in
$G$ between $A \setminus B$ and $B \setminus A$. The graph $G$ is
called $f$-connected if $|A \cap B| \geq f(|A \setminus B|)$, for
every separation $(A,B)$ of $G$ with $|A \setminus B| \leq |B
\setminus A|$. In~\cite{BBDK} it was proved that if $f(k) \geq
2(k+1)^2$ for every $k \in \mathbb{N}$ then $G$ is Hamiltonian for
every $n \geq 3$. It was also conjectured that there exists a
function $f$ which is linear in $k$ and is enough to ensure
Hamiltonicity. Using Theorem~\ref{th}, we can get quite close to
proving this conjecture for sufficiently large $n$:

\begin{theorem} \label{fconnected}
If $G=(V,E)$, where $|V|=n$, is $f$-connected for $f(k) = k \log k +
O(1)$, then it is Hamiltonian for sufficiently large $n$.
\end{theorem}

\noindent For the sake of simplicity and clarity of presentation, we
do not make a particular effort to optimize the constants obtained
in theorems we prove. We also omit floor and ceiling signs whenever
these are not crucial. All of our results are asymptotic in nature
and whenever necessary we assume that $n$ is sufficiently large.
Throughout the paper, $\log$ stands for the natural logarithm. We
say that some event holds \emph{almost surely}, or a.s. for brevity,
if the probability it holds tends to 1 as $n$ tends to infinity. Our
graph-theoretic notation is standard and follows that
of~\cite{Diestel}.

\noindent The rest of the paper is organized as follows: in
Section~\ref{sec::main} we prove and discuss Theorem~\ref{th}, in
Section~\ref{sec::cor} we prove its corollaries:
Theorems~\ref{uvpath},~\ref{pancyclicP2},~\ref{randomgraph}
and~\ref{fconnected}.

\section{Proof of the main result} \label{sec::main}

The proof of Theorem~\ref{th} is based on the ingenious
rotation-extension technique, developed by P\'osa~\cite{Posa}, and
applied later in a multitude of papers on Hamiltonicity (mostly of
random graphs). Our proof technique borrows some technical ideas
from the paper of Ajtai, Koml\'os and Szemer\'edi~\cite{AKS}.

Before diving into fine details of the proof, we would like to
compare our Hamiltonicity criterion and its proof with its
predecessors. Several previous papers, including~\cite{AKS},
~\cite{FK},~\cite{KrSu}, state, explicitly or implicitly, sufficient
conditions for Hamiltonicity applicable in principle to sparse
graphs. Usually criteria of this sort are carefully tailored to be
applied to random or pseudo-random graphs, and are therefore rather
complicated and not always natural. Moreover, such criteria are
sometimes fragile in the sense that they can be violated by adding
more edges to the graph -- a somewhat undesirable feature. Our
criterion in Theorem~\ref{th} is (on a qualitative level, at least)
quite natural and easily comprehensible, and can be potentially
applied to a very wide class of graphs.
As for our proof, due to the relative simplicity of the conditions we
use, the argument is perhaps more involved than some of the previous
proofs; there are however similarities. A novel
ingredient, relying heavily on Property P2, is the part presented in
Section~\ref{goodv} (finding many good initial rotations).

In order to be able to refer to the proof of our criterion while
proving some of the corollaries we break the proof into four parts,
each time indicating which property is needed for which part.

\begin{prop} \label{connected}
Let $G$ satisfy properties P1 and P2. Then $G$ is connected.
\end{prop}

\proofstart If not, let $C$ be the smallest connected component of
$G$. Then by P1, $|C|>\frac{n}{m}$, but then by P2, $E(C,V\setminus
C)\ne\emptyset$ -- a contradiction. \proofend

\subsection{Constructing an initial long path} \label{ILP}
In this subsection we show that a graph which satisfies some
expansion properties (that is, property P1 and some expansion of
larger sets, implied by property P2) contains a long path, and even
more, it has many paths of maximum length starting at the same
vertex.

Let $P_0=(v_1,v_2,\ldots,v_q)$ be a path of maximum length in $G$.
If $1 \leq i \leq q-2$ and $(v_q,v_i)$ is an edge of $G$ then
$P'=(v_1 v_2\ldots v_i v_q v_{q-1} \ldots v_{i+1})$ is also of
maximum length. $P'$ is called a {\em rotation} of $P_0$ with {\em
fixed endpoint} $v_1$ and {\em pivot} $v_i$. The edge
$(v_i,v_{i+1})$ is called the {\em broken} edge of the rotation. We
say that the segment $v_{i+1} \ldots v_q$ of $P_0$ is reversed in
$P'$.

In case the new endpoint, $v_{i+1}$, has a neighbor $v_j$ such that
$j \notin \{i, i+2\}$, then we can rotate $P'$ further to obtain
more paths of maximum length. We use rotations and extensions
together with property P1 to find a path of maximum length with
large rotation endpoint sets (see for
example~\cite{BFF1},~\cite{FK},~\cite{KS},~\cite{KrSu}).







\begin{claim}\label{St}
Let $G=(V,E)$ be a graph on $n$ vertices that satisfies property P1
and moreover any subset of $V$ of size $n/4130m$ has at least
$n-o(n)$ external neighbors. Let $P_0=(v_1,v_2,\ldots,v_q)$ be a path
of maximum length in $G$. Then there exists a set $B(v_1)\subseteq
V(P_0)$ of at least $n/3$ vertices, such that for every $v\in
B(v_1)$ there is a $v_1v$-path of maximum length which can be
obtained from $P_0$ by at most $\frac{2\log n}{\log d}$ rotations
with fixed endpoint $v_1$. In particular $|V(P_0)| \geq n/3$.
\end{claim}

\proofstart Let $t_0$ be the smallest integer such that
$\left(\frac{d}{3}\right)^{t_0-2} > \frac{n}{md}$. Note that $t_0
\leq  2\frac{\log n}{\log d}$.

We prove that there exists a sequence of sets $S_0, \ldots , S_{t_0}
= B(v_1) \subseteq V(P_0)\setminus\{v_1\}$ of vertices such that for
every $0 \leq t \leq t_0$, every $v\in S_t$ is the endpoint of a
path, obtainable from $P_0$ by $t$ rotations with fixed endpoint
$v_1$, such that for every $0 \leq i \leq t$, after the $i$th
rotation the non-$v_1$-endpoint of the path is in $S_i$, and
moreover $|S_t|=\left(\frac{d}{3}\right)^{t}$ for every $t\leq
t_0-3$, $|S_{t_0-2}|=\frac{n}{dm}$, $|S_{t_0-1}|=\frac{n}{4130m}$,
and $|S_{t_0}|\geq n/3$.

First we construct the sets by induction on $t$. For $t=0$, one can
choose $S_0=\{ v_q\}$ and all requirements are trivially satisfied.








Induction step: let $0< t\leq t_0-2$ and assume that the appropriate
sets $S_0, \ldots , S_{t-1}$ with the appropriate properties were
already constructed. We will now construct $S_t$. Let first

$$
T= \{v_i\in N(S_{t-1}) : v_{i-1},v_i,v_{i+1}\not\in
\bigcup_{j=0}^{t-1}S_j\}\ .
$$

be the set of potential pivots for the $t$th rotation. Assume now
that $v_i\in T$, $y\in S_{t-1}$ and $(v_i,y)\in E$. Then a $v_1
y$-path $Q$ can be obtained from $P_0$ by $t-1$ rotations such that
after the $j$th rotation, the non-$v_1$-endpoint is in $S_j$ for
every $j \leq t-1$. Each such rotation breaks an edge incident with
the new endpoint. Since $v_{i-1}, v_i, v_{i+1}$ are not endpoints
after any of these $t-1$ rotations, both edges $(v_{i-1},v_i)$ and
$(v_i,v_{i+1})$ of the original path $P_0$ must be unbroken and thus
must be present in $Q$.

Hence, rotating $Q$ with pivot $v_i$ will make either $v_{i-1}$, or
$v_{i+1}$ an endpoint (which one, depends on whether the unbroken
segment $v_{i-1} v_i v_{i+1}$ is reversed or not after the first
$t-1$ rotations). Assume w.l.o.g. it is $v_{i-1}$. We add $v_{i-1}$
to the set $\hat{S}_{t}$ of new endpoints and say that {\em $v_i$
placed $v_{i-1}$ in $\hat{S}_{t}$}. The only other vertex that can
place $v_{i-1}$ in $\hat{S}_{t}$ is $v_{i-2}$ (if it exists). Thus,

\begin{eqnarray*}
|\hat{S}_{t}| &\ge& \frac{1}{2} |T|\ge
\frac{1}{2} \left(|N(S_{t-1})|-3(1+|S_1|+\ldots+|S_{t-1}|\right))\\
&\geq& \frac{d}{2} \left(\frac{d}{3}\right)^{t-1} -
\frac{3}{2}\frac{(d/3)^{t}-1}{d/3-1} \geq
\left(\frac{d}{3}\right)^{t}
\end{eqnarray*}

where the last inequality follows since $d \geq 12$. Clearly we can
delete arbitrary elements of $\hat{S}_{t}$ to obtain $S_{t}$ of size
exactly $\left(\frac{d}{3}\right)^{t}$ if $t\leq t_0-3$ and of size
exactly $\frac{n}{dm}$ if $t=t_0-2$. So the proof of the induction
step is complete and we have constructed the sets $S_0, \ldots,
S_{t_0-2}$.

To construct $S_{t_0-1}$ and $S_{t_0}$ we use the same technique as
above, only the calculation is slightly different. Since
$|N(S_{t_0-2})|\geq d\cdot \frac{n}{dm}$, we have

\begin{eqnarray*}
|\hat{S}_{t_0-1}| &\ge& \frac{1}{2} |T|\ge
\frac{1}{2} \left(|N(S_{t_0-2})|-3(1+|S_1|+\ldots+|S_{t_0-3}|+|S_{t_0-2}|\right))\\
&\geq& \frac{n}{2m} -
\frac{3}{2}\left(\frac{(d/3)^{t_0-3}-1}{(d/3)-1}
+2\frac{n}{dm}\right) \geq
\frac{n}{2m} - \frac{3}{2}\cdot\left(\frac{d}{3}\right)^{t_0-3} -3\frac{n}{dm}\\
&\geq& \frac{n}{2m} - \frac{3}{2}\cdot\frac{n}{dm} -3\frac{n}{dm}
\geq \frac{n}{4130m},
\end{eqnarray*}

where the last inequality follows since $d \geq 12$.

For $S_{t_0}$ the difference in the calculation comes from using the
expansion guaranteed by property P2 rather than the one guaranteed
by property P1, that is, $|N(S_{t_0-1})| \geq n-o(n)$. We have

\begin{eqnarray*}
|{S}_{t_0}| &\ge& \frac{1}{2} |T|\ge
\frac{1}{2} \left(|N(S_{t_0-1})|-3(1+|S_1|+\ldots+|S_{t_0-2}|+|S_{t_0-1}|\right))\\
&\geq& \frac{n}{2}(1-o(1)) -
\frac{3}{2}\left(\frac{(d/3)^{t_0-3}-1}{(d/3) -1}
+\frac{2n}{d m}+\frac{n}{4130m}\right)\\
&\geq& \frac{n}{2}(1-o(1)) -
\frac{3}{2}\left(\frac{3n}{d m}+\frac{n}{4130m}\right)\\
&\geq& \frac{n}{3},
\end{eqnarray*}

where the last inequality follows since $n\geq 2m$ and $d\geq 12$.

The set $S_{t_0}$ can be chosen to be $B(v_1)$ and satisfies all the
requirements of the Claim. Note that since $S_{t_0}\subseteq
V(P_0)$, we have $|V(P_0)|> n/3$. This concludes the proof of the
claim. \proofend

{\bf Remark} Note that, although we do not need it here, the
rotations which create these paths always brake an edge of the
original path $P_0$.

\subsection{Finding many good initial rotations}\label{goodv}

In this subsection we prove an auxiliary lemma, which will be used
in the next subsection to conclude the proof of Theorem~\ref{th}.

Let $H$ be a graph with a spanning path $P=(v_1,\ldots,v_l)$. For
$2\leq i < l$ let us define the auxiliary graph $H^+_i$ by adding a
vertex and two edges to $H$ as follows: $V(H^+_i)= V(H)\cup \{ w\}$,
$E(H^+_i)=E(H)\cup \{(v_l,w), (v_i,w)\}$. Let $P_i$ be the spanning
path of $H^+_i$ which we obtain from the path $P\cup\{(v_l,w)\}$ by
rotating with pivot $v_i$. Note that the endpoints of $P_i$ are
$v_1$ and $v_{i+1}$.


For a vertex $v_i \in V(H)$ let $S^{v_i}$ be the set of those
vertices of $V(P)\setminus \{ v_1\}$, which are endpoints of a
spanning path of $H_i^+$ obtained from $P_i$ by a series of
rotations with fixed endpoint $v_1$.

A vertex $v_i\in V(P)$ is called a {\em bad initial pivot} (or
simply a {\em bad vertex}) if $|S^{v_i}|< \frac{l}{43}$ and is
called a {\em good initial pivot} (or a {\em good vertex})
otherwise. We can rotate $P_i$ and find a large number of endpoints
provided $v_i$ is a good initial pivot.


Using an argument similar to the one used in the proof of
Claim~\ref{St}, we can show that $H$ has many good initial pivots
provided that a certain condition, similar to property P2, is
satisfied.

\begin{lemma}\label{good}
Let $H$ be a graph with a spanning path $P=(v_1,\ldots,v_l)$. Assume
that every two disjoint sets $A$, $B$ of vertices of $H$ of sizes
$|A|, |B|\ge l/43$ are connected by an edge. Then we have $$|R| \le
7l/43,$$ where $R=R(P)\subseteq V(P)$ is the set of bad vertices.
\end{lemma}

\proofstart We will create a set $U\subseteq V(H)$, whose size is at
least $|R|/7$, but does not expand enough, that is, $|U\cup
N_H(U)|\leq 21|U|$. This in turn will imply that the set $R$ of bad
vertices cannot be big.

Let $R=\{v_{i_1},\ldots, v_{i_r}\}$. We {\em process} the vertices
of $R$ one after the other. We will maintain subsets $U$ and $X$ of
$V(H)$ where initially $U=X=\emptyset$. Whenever we finish
processing a vertex of $R$ we update the sets $U$ and $X$. The
following properties will hold after the processing of $v_{i_j}$.

\begin{equation}\label{UX}
U \subseteq X \,,\quad N_H(U)\subseteq ext(X)\,,\quad |U|\ge
\frac{1}{7}|X| \,,\quad \{v_{i_1+1},\ldots,v_{i_j+1}\}\subseteq X \
,
\end{equation}

where $ext(X)$ denotes the set containing the vertices of $X$
together with their left and right neighbors on $P$. Clearly
$|ext(X)|\leq 3|X|$.

Suppose the current vertex to process is $v_{i_j}$. If $v_{i_j+1}\in
X$, then we do not change $U$ and $X$ and so the conditions of
(\ref{UX}) trivially hold by induction.

Otherwise, we will create sets $W_t\subseteq S^{v_{i_j}}$
inductively, such that for every $t$ the following hold.
\begin{itemize}


\item[$(a)$] $ W_t\subseteq S_t^{v_{i_j}}$;\\[-5mm]

\item[$(b)$] $ |W_t|=2^t $;\\[-5mm]

\item[$(c)$] $W_t\cap \left(\cup_{s=0}^{t-1}W_s\cup X\right)=\emptyset$,

\end{itemize}

where $S_t^{v_{i_j}}$ contains those vertices $y$ of $S^{v_{i_j}}$
for which a spanning path of $H^+_i$ ending at $y$ can be produced
from $P_i$ by $t$ rotations with fixed endpoint $v_1$, such that
after the $s$th rotation the new endpoint is in $W_s$, for every $s
< t$.

We begin by setting $W_0=\{v_{i_j+1}\}$. Conditions $(a)$ and $(b)$
trivially hold, for condition $(c)$ note that $v_{i_j+1}\notin X$.

Assume now that we have constructed $W_0, \ldots , W_t$ with
properties $(a)-(c)$. If $|N_H(W_t)\setminus
ext(\left(\cup_{i=1}^{t}W_i\cup X\right))| > 5|W_t|$, then we create
$W_{t+1}$ with properties $(a)-(c)$, otherwise we finish the
processing of $v_{i_j}$ by updating $U$ and $X$.

Let $T_t=N_H(W_t)\setminus ext(\left(\cup_{i=1}^{t}W_i\cup
X\right))$ and assume first that $|T_t| >  5|W_t|$. We use an
argument similar to the one used in Claim~\ref{St} to create
$W_{t+1}$ with properties $(a)-(c)$.









Let $v_i\in T_t$, $v_i\neq v_1,v_l$, and suppose that $v_i$ is
adjacent to $y\in W_t$. Recall, that by property $(a)$ a spanning
path $Q$ of $H^+_i$ ending at $y$ can be produced from $P_i$ by $t$
rotations, such that for every $s < t$, after the $s$th rotation the
new endpoint is in $W_s$. Since the vertices $v_{i-1}, v_i$ and
$v_{i+1}\not\in\bigcup_{s=0}^t W_s$, they are not endpoints after
any of these $t$ rotations. Each rotation breaks an edge incident
with the new endpoint, hence both edges $(v_{i-1},v_i)$ and
$(v_i,v_{i+1})$ of the original path $P_i$ must be  present in $Q$.
Rotating $Q$ with pivot $v_i$ will brake one of them. Such a
rotation also makes one of $v_{i-1}$ and $v_{i+1}$ into an endpoint,
and as such, into an element of $S^{v_{i_j}}_{t+1}$. Denote this
vertex by $v_{i}'$. We define $W_{t+1}=\{v_i' : v_i\in  T_t\}$. We
say that {\em $v_i'$ is placed in $W_{t+1}$ by $v_i$}. Observe that
besides $v_i$ the only other vertex that can place $v_{i}'$ in
$W_{t+1}$ is its other neighbor on the path $P_i$. Thus,
$$
|W_{t+1}| \ge \left\lceil\frac{1}{2} (|T_t|-2)\right\rceil \geq
2|W_t|.
$$

Deleting arbitrarily some vertices from $W_{t+1}$ we can make sure
that its cardinality is exactly $2|W_t|$. Properties $(a)$ and $(b)$
are then naturally satisfied. Property $(c)$ is satisfied because,
by the definition of $T_t$ we have $v_i\notin
ext\left(\cup_{s=0}^{t}W_s\cup X\right)$ and so none of its
neighbors on $P_i$, in particular $v_i'$, is an element of
$\left(\cup_{s=0}^{t}W_s\cup X\right)$.

Property $(b)$ ensures that $|\cup_{s=0}^{t} W_s|$ is strictly
increasing so the processing of the vertex $v_{i_j}$ is bound to
reach a point in which $|T_k| \leq 5|W_k|$ for some index $k$. At
that point we update $U$ and $X$ by adding $W_k$ to $U$ and adding
$W_1\cup \cdots \cup W_k\cup T_k$ to $X$. We have to check that the
conditions of (\ref{UX}) hold.

Observe that $|W_1\cup \cdots \cup W_k| <2|W_k|$, so the number of
vertices added to $X$ is at most seven times more than the number of
vertices added to $U$. Also, property $(c)$ and $U\subseteq X$ made
sure that $W_k$ was disjoint from $U$, so indeed the property
$|U|\geq |X|/7$ remains valid. The other conditions in (\ref{UX})
follow easily from the definition of the ``new'' $U$ and $X$. Hence
the processing of $v_{i_j}$ is complete.

{\bf Claim} $|U|\leq l/43$.

\proofstart Assume the contrary and let $j$ be the smallest index,
such that $|U|> l/43$ after the processing of $v_{i_j}$.\\
Observe that $|U|\le 2l/43$. Indeed, after the processing of
$v_{i_j}$ the set $U$ received at most $|S^{v_{i_j}}|$ vertices,
which is at most $l/43$, due to the fact that $v_{i_j}$ is a bad
vertex. We thus have $l/43<|U|\le 2l/43$, $U\subseteq X$,
$N_H(U)\subseteq ext(X)$ and $|ext(X)|\leq 3|X|\leq 21|U|$. Then
$|V(P)\setminus ext(X)|\ge l/43$, and there are no edges of $H$
between $U$ and $V(P)\setminus ext(X)$. This contradicts our
assumption on $H$. \proofend

To conclude the proof of the Lemma we note that after processing all
vertices of $R$, we have
$R^+:=\{v_{i_1+1},\ldots,v_{i_r+1}\}\subseteq X$ and $|U|\ge |X|/7$
by (\ref{UX}). Since $|U|\le l/43$, it follows that $|R|=|R^+|\le
7l/43$. \proofend

\subsection{Closing the maximal path} \label{sec::close}

\begin{lemma}
Let $G$ be a connected graph that satisfies property P2. Let the
conclusion of Claim~\ref{St} be also true for $G$, that is, for
every path $P_0=(v_1,v_2,\ldots,v_q)$ of maximum length in $G$ there
exists a set $B(v_1)\subseteq V(P_0)$ of at least $n/3$ vertices,
such that for every $v\in B(v_1)$ there is a $v_1v$-path of maximum
length which can be obtained from $P_0$ by at most $t_0 \leq
\frac{2\log n}{\log d}$ rotations with fixed endpoint $v_1$. Then
$G$ is Hamiltonian.
\end{lemma}

\proofstart We will prove that there exists a path of maximum length
which can be closed into a cycle. This, together with connectedness
implies that the cycle is Hamiltonian. To find such a path of
maximum length we will create two sets of vertices, large enough to
satisfy property P2, such that between any two vertices (one from
each) there is a path of maximum length.

Let $P_0=(v_1,v_2,\ldots,v_q)$ be a path of maximum length in $G$.
Let $A_0=B(v_1)$. For every $v\in A_0$ fix a $v_1v$-path $P^{(v)}$
of maximum length and, using our assumption, construct sets $B(v)$,
$|B(v)| \geq n/3$, of endpoints of maximum length paths with
endpoint $v$, obtained from a $P^{(v)}$ by at most $t_0$ rotations.
In summary, for every $a\in A_0$, $b\in B(a)$ there is a maximum
length path $P(a,b)$ joining $a$ and $b$, which is obtainable from
$P_0$ by at most $\rho:=2t_0 \leq \frac{4\log n}{\log d}$ rotations.

We consider $P_0$ to be directed and divided into $2\r$ segments
$I_1,I_2,\ldots,I_{2\r}$ of length at least $\rdown{|P_0|/2\r}$
each, where $|P_0|\ge n/3$. As each $P(a,b)$ is obtained from $P_0$
by at most $\r$ rotations and every rotation breaks exactly one edge
of $P_0$, the number of segments of $P_0$ which occur complete on
this path, although perhaps reversed, is at least $\r$. We say that
such a segment is {\em unbroken}. These segments have an absolute
orientation given to them by $P_0$, and another, relative to this
one, given to them by $P(a,b)$, which we consider to be directed
from $a$ to $b$. We consider sequences $\s=I_{i_1}, I_{i_2}, \ldots,
I_{i_\tau}$ of unbroken segments of $P_0$, which occur in this order
on $P(a,b)$, where $\s$ also specifies the relative orientation of
each segment. We call such a sequence $\s$ a {\em $\tau$-sequence},
and say that $P(a,b)$ {\em contains} $\s$.

\noindent For a given $\tau$-sequence $\s$, we consider the set
$L(\s)$ of ordered pairs $(a,b)$, $a\in A_0, \; b\in B(a)$, such
that $P(a,b)$ contains $\s$.

The total number of $\tau$-sequences is $2^\tau (2\r)_\tau$. Any
path $P(a,b)$ contains at least $\r$ unbroken segments, and thus at
least $\binom{\r}{\tau}$ $\tau$-sequences. The average, over
$\tau$-sequences, of the number of pairs $(a,b)$ such that $P(a,b)$
contains a given $\tau$-sequence is therefore at least
$$\frac{n^2}{9} \cdot \frac{\binom{\r}{\tau}}{2^\tau (2\r)_\tau} \ge \a n^2,$$ where
$\a = \a(\tau) = 1/9(4\tau)^{-\tau}$. Thus, there is a
$\tau$-sequence $\s_0$ and a set $L= L(\s_0),\,|L|\geq \a n^2$ of
pairs $(a,b)$ such that for each $(a,b)\in L$ the path $P(a,b)$
contains $\s_0$. Let $\hat{A}=\{a \in P_0 : L$ contains at least $\a
n/2$ pairs with $a$ as first element\}. Then $|\hat{A}|\geq \a n/2$.
For each $a\in \hat{A}$ let $\hat{B}(a)=\{b : (a,b)\in L\}$. Then,
by the definition of $\hat{A}$, for each $a\in\hat{A}$ we have
$|\hat{B}(a)|\ge \a n/2$.

Let $\tau = \frac{\log \log n}{2\log \log \log n}$ and let
$\s_0=(I_{i_1}, I_{i_2}, \ldots, I_{i_\tau})$. We divide $\sigma_0$
into two sub-sequences, $\s_0^1 = (I_{i_1}, \ldots, I_{i_{\tau/2}})$
and $\s_0^2 = (I_{i_{\tau/2+1}}, \ldots, I_{i_\tau})$ where both
sub-sequences maintain the order and orientation of the segments of
$\s_0$. Both sub-sequences $\s_0^1$ and $\s_0^2$ have at least
$\tau/2 \cdot n/(6\r)\ge \frac{n}{96m}$ vertices. Let $x$ be the
last vertex of $I_{i_{\tau/2}}$, and let $y$ be the first vertex of
$I_{i_{\tau/2+1}}$ (in the orientation given by $\s_0$). Now we
define the notion of good vertices in $\s_0^1$ and $\s_0^2$. For
$\s_0^1$ construct a graph $H_1$ from the segments of $\sigma_0^1$
by joining by an edge the last vertex of $I_{i_j}$ to the first
vertex of $I_{i_{j+1}}$ for every $1 \leq j < \tau/2$ and then
adding the edges of $G$ with both endpoints in the interior (that
is, not endpoints) of segments of $\s_0^1$ to $H_1$. Then the
segments of $\s_0^1$ with the edges linking them form an oriented
spanning path in $H_1$, starting at $x$. We define {\em good}
vertices in $\s_0^1$ to be the vertices which are not endpoints of
any segment of $\s_0^1$ and are good vertices of $H_1$ as defined in
Section \ref{goodv}, with $l=\sum_{j=1}^{\tau/2}|I_{i_j}|$. Due to
property P2, Lemma \ref{good} applies here, and so, since $\tau =
o(|\s_0^1|)$, more than half of the vertices of $\s_0^1$ are good.
For $\s_0^2$ we act similarly: construct a graph $H_2$ from the
segments of $\sigma_0^2$ by joining the first vertex of $I_{i_j}$ to
the last vertex of $I_{i_{j-1}}$ for every $\tau/2+1 < j \leq \tau$
and then adding the edges of $G$ with both endpoints in the interior
of segments of $\s_0^2$ to $H_2$. Then the segments of $\s_0^2$ with
the edges linking them and $h_2$ form an oriented spanning path in
$H_2$, starting at $y$. We define {\em good} vertices in $\s_0^2$ to
be the vertices which are not endpoints of any segment of $\s_0^2$
and are good vertices of $H_2$ as defined in Section \ref{goodv},
with $l=\sum_{j=\tau/2+1}^{\tau}|I_{i_j}|$. Due to property P2,
Lemma \ref{good} applies here, and so, since $\tau = o(|\s_0^2|)$,
more than half of the vertices of $\s_0^2$ are good.

Since $|\hat{A}| \geq \a n/2 \geq \frac{n}{4130 m}$ (which is why we
get the upper bound on $d$ in Theorem~\ref{th}) and $\s_0^1$ has at
least $|\s_0^1|/2 > \frac{n}{192m}$ good vertices, there is an edge
from a vertex $\hat{a}\in \hat{A}$ to a good vertex in $\s_0^1$.
Similarly, as $|\hat{B}(\hat{a})|\ge \a n/2$ ,there is an edge from
some $\hat{b}\in \hat{B}(\hat{a})$ to a good vertex in $\s_0^2$.
Consider the path $\hat{P}=P(\hat{a},\hat{b})$ of maximum length
connecting $\hat{a}$ and $\hat{b}$ and containing $\sigma_0$. The
vertices $x$ and $y$ split this path into three sub-paths: $P_1$
from $\hat{a}$ to $x$, $P_2$ from $y$ to $\hat{b}$ and $P_3$ from
$x$ to $y$. We will rotate $P_1$ with $x$ as a fixed endpoint and
$P_2$ with $y$ as a fixed endpoint. We will show that the obtained
endpoint sets $V_1$ and $V_2$ are sufficiently large. Then by
property P2 there will be an edge of $G$ between $V_1$ and $V_2$.
Since we did not touch $P_3$, this edge closes a maximum path into a
cycle, which is Hamiltonian due to the connectivity of $G$.

Since there is an edge from $\hat{a}$ to a good vertex in $\s_0^1$,
by the definition of a good vertex we can rotate $P_1$, starting
from this edge, to get a set $V_1$ of at least $|\s_0^1|/43 >
n/(4130m)$ endpoints. When doing this, we will treat the subpath
that links $\hat{a}$ and the first vertex of $I_{i_1}$ and each
subpath that links two consecutive segments of $\sigma_0^1$, as
single edges and ignore edges of $G$ that are incident with an
endpoint of some segment of $\sigma_0^1$ - like in $H_1$. This
ensures that all rotations and broken edges are inside segments of
$\sigma_0^1$ and so there is indeed a path of the appropriate length
from $x$ to every vertex of $V_1$.

Similarly, since there is an edge from $\hat{b}$ to a good vertex in
$\s_0^2$, we can rotate $P_2$, starting from this edge to get a set
$V_2$ of at least $|\s_0^2|/43 > n/(4130m)$ endpoints. When doing
this, we will treat the subpath that links $\hat{b}$ and the last
vertex of $I_{i_{\tau}}$ and each subpath that links two consecutive
segments of $\sigma_0^2$, as single edges and ignore edges of $G$
that are incident with an endpoint of some segment of $\sigma_0^2$ -
like in $H_2$. This ensures that all rotations and broken edges are
inside segments of $\sigma_0^2$ and so there is indeed a path of the
appropriate length from $y$ to every vertex of $V_2$. This concludes
the proof of Theorem~\ref{th}. \proofend

\subsection{Hamiltonicity with larger expansion}\label{large-d}
As we have mentioned, our Hamiltonicity criterion can be extended to
handle graphs with a larger expansion than that postulated in
Theorem \ref{th} ($d\le e^{\sqrt[3]{\log n}}$). In particular, using
very similar arguments, we can prove the following statement.

\begin{theorem}\label{th2}
Let $12 \leq d \leq \sqrt{n}$ and let $G$ be a graph on $n$ vertices
satisfying the following two properties:
\begin{description}
\item[P1'] For every $S \subset V$, if $|S|\le \frac{n \log
d}{d\log n}$ then $|N(S)|\ge d|S|$;
\item[P2'] There is an edge in $G$ between any two disjoint
subsets $A,B \subseteq V$ such that $|A|,|B|\ge \frac{n \log d}{1035
\log n}$.
\end{description}

Then $G$ is Hamiltonian, for sufficiently large $n$.
\end{theorem}

The proof of Theorem \ref{th2} is almost identical to that of
Theorem \ref{th} given above. The only notable difference is that
here we can allow ourselves to take $\tau=2$ in the proof.

\section{Corollaries} \label{sec::cor}

In this section we prove the afore-mentioned corollaries of
Theorem~\ref{th}.

\textbf{Proof of Theorem~\ref{uvpath}} Let $G_{uv} = (V,E \cup
\{(u,v)\})$; clearly $G_{uv}$ satisfies properties P1 and P2 and is
therefore Hamiltonian by Theorem~\ref{th}. Let $C = w_1 w_2 \ldots
w_n w_1$ be a Hamilton cycle in $G_{uv}$. If $(u,v)$ is an edge of
$C$, remove it to obtain the desired path in $G$. Otherwise,
assuming that $u = w_i$ and $v = w_j$, add $(u,v)$ to $E(C)$ and
remove $(u,w_{i+1})$ and $(v,w_{j+1})$, where all indices are taken
modulo $n$, to obtain a Hamilton path of $G_{uv}$ that contains the
edge $(u,v)$; denote this path by $P$. We will close $P$ into a
Hamilton cycle that includes $(u,v)$; removing this edge will result
in the required path. The building of the cycle will be done as in
the proof of Theorem~\ref{th} Section~\ref{sec::close}, with $P$ as
$P_0$, while making sure that $(u,v)$ is never broken. The proof is
essentially the same, except for the following minor changes:

\begin{enumerate}
\item When dividing $P$ into $2\r$ segments, we will make sure that
$(u,v)$ is in one of the segments; denote it by $I_j$.

\item When considering $\tau$-sequences, we will restrict ourselves
to those that include $I_j$.

\item Assume without loss of generality that $I_j \in \sigma_0^1$.
When building $H_1$ (and later, when rotating $P_1$ according to the
model of $H_1$) we will ignore $I_j$, that is, we will replace it by
a single edge $(a,b)$ where $a$ is the last vertex of $I_{j-1}$ (or
$h_1$ if $j=1$) and $b$ is the first vertex of $I_{j+1}$ (or $x$ if
$j = \tau/2$).
\end{enumerate}

{\hfill $\Box$
\medskip\\}

\textbf{Proof of Theorem~\ref{pancyclicP2}}

Fix some $\frac{8 n \log n}{t \log \log n} \leq k \leq n - 3n/t$.
Let $V_0 \subseteq V$ be an arbitrary subset of size $k + n/t$. We
construct a sequence of subsets $S_i$, let $S_0 = \emptyset$. As
long as $|S_i| < n/t$ and there exists a set $A_i \subseteq V_0
\setminus S_i$ such that $|A_i| \leq n/t$ but $|N_{G[V_0 \setminus
S_i]}(A_i)| < |A_i|\frac{4\log n}{\log \log n}$, we define $S_{i+1}
:= S_i \cup A_i$. Let $q$ be the smallest integer such that $|S_q|
\geq n/t$ or $|N_{G[V_0 \setminus S_q]}(A)| \geq |A|\frac{4\log
n}{\log \log n}$ for every $A \subseteq V_0 \setminus S_q$ of size
at most $n/t$. We claim that $|S_q| < n/t$. Indeed assume for the
sake of contradiction that $|S_q| \geq n/t$. Since we halt the
process as soon as this occurs, and $|A_{q-1}| \leq n/t$, we have
$|S_q| < 2 n/t$. For every $0 \leq i \leq q-1$ we have $|N_{G[V_0
\setminus S_i]}(A_i)| < |A_i|\frac{4\log n}{\log \log n}$ and so
$|N_{G[V_0]}(S_q)| < |S_q|\frac{4\log n}{\log \log n}$. On the other
hand, $G$ satisfying property P$2$ together with our lower bound on
$k$ implies $|N_{G[V_0]}(S_q)|
> |V_0| - n/t - |S_q| \geq |V_0| - 3n/t \geq k \geq |S_q|\frac{4\log n}{
\log \log n}$, a contradiction.

Hence, $|S_q| < n/t$ and so, for $U = V_0 \setminus S_q$, $G[U]$
satisfies an expansion condition similar to {\bf P1}, that is, for
every $A \subseteq U$, if $|A| \leq n/t$ then $|N_{G[U]}(A)| \geq
4|A|\frac{\log n}{\log \log n}$.

In the following we prove that with positive probability the induced
subgraph of $G$ on a random $k$-element subset of $U$ also satisfies
a condition similar to {\bf P1}. Let $K$ be a $k$-subset of $U$
drawn uniformly at random. We will prove that, with positive
probability, $G[K]$ satisfies the following:
\begin{description}
\item[P1] For every $A \subseteq K$, if $|A| \leq n/t$
then $|N_{G[K]}(A)| \geq 2|A|\frac{\log n}{\log \log n}$.
\end{description}

Let $r = |U|-k$. Note that $0 \leq r \leq n/t$. Let $A \subseteq U$
be any set of size $a \leq n/t$, then, as was noted above,
$|N_{G[U]}(A)| \geq 4|A|\frac{\log n}{\log \log n}$. Let $N_0
\subseteq N_{G[U]}(A)$ be an arbitrary subset of size
$4|A|\frac{\log n}{\log \log n}$. If $A \subseteq K$ and
$|N_{G[K]}(A)| \leq 2|A|\frac{\log n}{\log \log n}$, then $K$ misses
at least $2|A|\frac{\log n}{\log \log n}$ vertices from $N_0$. This
can occur with probability at most

\begin{eqnarray*}
\frac{{|N_0| \choose \frac{2a \log n}{\log \log n}} {|U| - \frac{2a
\log n}{\log \log n} \choose r - \frac{2a \log n}{\log \log
n}}}{{|U| \choose r}} &\leq& {\frac{4a \log n}{\log \log n} \choose
\frac{2a \log n}{\log \log n}}\left(\frac{r}{|U|}\right)
^{\frac{2a \log n}{\log \log n}}\\
&\leq& 2^{\frac{4a \log n}{\log \log n}}\left(\frac{\frac{
n}{t}}{\frac{8 n \log n}{t \log \log n}}\right)^{\frac{2a \log n}{\log \log n}}\\
&=& \left(\frac{\log \log n}{2 \log n}\right)^{\frac{2a \log n}{\log
\log n}}.
\end{eqnarray*}
Note that the latter bound is $o(\frac{1}{n})$ for $a=1$, and
$o(\frac{1}{n} {n \choose a}^{-1})$ for every $a \geq 2$.

It follows by a union bound argument that
$$
Pr\left[\text{ there exists an } A \subseteq K \text{ such that }
|A| \leq n/t \text{ but } N_{G[K]}(A) < \frac{2\log n}{\log \log
n}|A|\right] = o(1).
$$
Hence, there exists an $k$-subset $X$ of $U$ such that for every $A
\subseteq X$, if $|A| \leq n/t$ then $|N_{G[X]}(A)| \geq \frac{2\log
n}{\log \log n}|A|$. Moreover, if $A,B$ are disjoint subsets of $V$,
and $|A|,|B| \geq \frac{k \log \log k \log \left(\frac{2\log n}{\log
\log n}\right)}{4130 \log k \log \log \log k} \geq n/t$ then there
is an edge between a vertex of $A$ and a vertex of $B$.

Thus $G[X]$ satisfies the conditions of Theorem~\ref{th} with
$|V|=k$ and $d = \frac{2\log n}{\log \log n}$ and is therefore
Hamiltonian. It follows that $G$ admits a cycle of length exactly
$k$.

{\hfill $\Box$
\medskip\\}

\textbf{Proof of Theorem~\ref{randomgraph}} Let $G=G(n,p)=(V,E)$ and
let $d = (\log n)^{0.1}$. We begin by showing that a.s. $G$
satisfies property P2 with respect to $d$. Indeed

\begin{eqnarray*}
Pr[G \nvDash P2] &\leq& {n \choose \frac{n \log \log n \log d}{4130
\log n \log \log \log n}}^2 \left(1-\frac{\log n + \log \log n +
\omega(1)}{n}\right)^{\left(\frac{n \log \log n \log d}{4130 \log n
\log \log \log n}\right)^2}\\ &\leq& \left(\frac{4130 e \log n \log
\log \log n}{0.1 (\log \log n)^2}\right)^{\frac{0.2 n (\log \log
n)^2} {4130 \log n \log \log \log n}}\\ & & \times \exp \left\{-
\frac{\log n + \log \log n + \omega(1)}{n}
\cdot \frac{0.01 n^2 (\log \log n)^4}{4130^2 (\log n)^2 (\log \log \log n)^2}\right\}\\
&=& o(1).
\end{eqnarray*}

Next, we deal with property P1. Since a.s. there are vertices of
"low" degree in $G$, we cannot expect every "small" set to expand by
a factor of $d$. Therefore, to handle this difficulty, we introduce
some minor changes to the proof of Theorem~\ref{th}, in fact only to the part
included in Claim~\ref{St}. First of all,
note that a.s. $G$ is connected (this fact replaces
Proposition~\ref{connected}). Let $SMALL=\{u \in V : d_G(u) \leq
(\log n)^{0.2}\}$ denote the set of all vertices of $G$ that have a
"low" degree. The vertices in $SMALL$ will be called \emph{small
vertices}. Standard calculations show that a.s. $G$ satisfies the
following properties:

\begin{itemize}

\item[$(1)$] $\delta(G) \geq 2$.

\item [$(2)$] For every $u \neq v \in SMALL$ we have $dist_G(u,v) \geq
250$, where $dist_G(u,v)$ is the number of edges in a shortest path
between $u$ and $v$ in $G$.

\item [$(3)$] $G$ satisfies a weak version of P1, that is, if $A \subseteq
V \setminus SMALL$ and $|A| \leq \frac{n \log \log n \log d}{d\log n
\log \log \log n}$ then $|N_G(A)| \geq 3d|A|$.

\item [$(4)$] The number of vertices of degree at most 11 is $O(\log^{11} n)$.

\end{itemize}

We will prove that, based on these properties, we can build initial
long paths as in Claim~\ref{St} of the proof
of Theorem~\ref{th}; this will conclude our proof of
Theorem~\ref{randomgraph}, as in Subsections~\ref{goodv}
and~\ref{sec::close} we did not rely on property P1. The argument is
essentially the same as in Claim~\ref{St}; the main difference is
that we will use roughly twice as many rotations to create the
eventual endpoint set of size $n/3$. This extra factor two has no
real effect on the rest of the proof.

Suppose first that the initial path of maximum length $P_0$ is such
that, while creating the sets $S_1, \ldots , S_{120}$ as we did in
the proof of Claim~\ref{St}, no vertex from $\cup_{i=1}^{119}S_i$ is
a small vertex. Then, by $(3)$, like in the proof of Claim~\ref{St},
after the $i$th rotation there are exactly $(3d/3)^{i} = (\log
n)^{0.1i}$ new endpoints in $S_i$. Therefore, after $120$ rotations
we will have an endpoint set $S_{120}$ with $(\log n)^{12}$
elements.

Suppose now that there is a vertex $u\in S_j\cap SMALL$ for some
$j\leq 119$. Let $P_u$ denote a path of maximum length from $v_1$ to
$u$ (which can be obtained from $P_0$ by at most 119 rotations). At
this point we ignore the endpoint sets $S_i$, $i\leq j$ created so
far and restart creating them. The first rotation is somewhat
special. By property $(1)$, $u$ has at least one neighbor on $P_u$
other than its predecessor. Thus we can rotate $P_u$ once and obtain
a $v_1w$-path $P_w$ of maximum length, such that $w$ is at distance
two from a small vertex. We create new endpoint sets $S_1, \ldots,
S_{120}$ with $P_w$ as the initial path. Note that property $(2)$
implies $w \notin SMALL$. Since a new endpoint is always at distance
at most two from the old endpoint, we can rotate another 120 times
without ever creating an endpoint which is a small vertex. Thus,
property $(3)$ applies and after the $i$th rotation (not including
the one that turned $w$ into an endpoint), $i\leq 120$, there are
exactly $(3d/3)^{i} = (\log n)^{0.1i}$ new endpoints in $S_i$.
Hence, after 120 further rotations we obtain a set $S_{121}$ of size
exactly $(\log n)^{12}$. Altogether we used up to $240$ rotations.

In the following we will prove that the endpoint sets we build grow
by the same multiplicative factor every at most {\em two} rotations.

We will prove by induction on $t$ that there exist endpoint sets
$S_{121}, S_{122}, \ldots$ such that for every $t\geq 122$, either
$|S_t|=\frac{d}{3}|S_{t-1}|$ or
$|S_t|=|S_{t-1}|=\frac{d}{3}|S_{t-2}|$.

Note that this implies $\sum_{i=0}^{t} |S_i| \leq \frac{4}{3}|S_t|$
if $|S_t|=\frac{d}{3}|S_{t-1}|$, provided $n$ is large enough.

For the base case we just have to note that $\sum_{i=0}^{121} |S_i|
\leq \frac{4}{3}|S_{121}|$. Suppose we have already built $S_t$ for
some $t \geq 121$ such that $\sum_{i=0}^{t} |S_i| \leq
\frac{4}{3}|S_t|$ and now wish to build $S_{t+1}$. We will proceed
as in the proof of Claim~\ref{St}.

Assume first that $|N(S_t)| \geq  d|S_t|$. Then, as in the proof of
Claim~\ref{St}
$$|\hat{S}_{t+1}| \geq \frac{1}{2}(d |S_{t}|-
3\cdot \frac{4}{3}|S_t|) =\frac{d-4}{2}|S_t|.$$ Hence, a subset
$S_{t+1}\subseteq \hat{S}_{t+1}$ with $|S_{t+1}|=\frac{d}{3}|S_t|$
can be selected.

Assume now that $|N(S_t)| < d|S_t|$. By $(3)$, this must mean that
for $S_t' := S_t \cap SMALL$ we have $|S_t'| \geq \frac{2}{3}|S_t|$.
Since $|S_t'| \gg \log^{11} n$, property $(4)$ implies that almost
every vertex of $S_t'$ has degree at least 12. By $(3)$, no two
small vertices have a common neighbor, so $|N(S_t')| \geq
(12-o(1))|S_t'|\geq (8-o(1))|S_t|$. As in the proof of
Claim~\ref{St}, we have $$|\hat{S}_{t+1}|\geq \frac{1}{2}(|N(S_t')|
-3\cdot |\cup_{i=1}^t S_i|)\geq \frac{1}{2} ((8-o(1))|S_t|-3\cdot
\frac{4}{3}|S_t|) \geq |S_t|.$$ Hence we can select an
$S_{t+1}\subseteq \hat{S}_{t+1}$ such that $|S_{t+1}|=|S_t|$.
Crucially, since we only used vertices from $S_t'$ for further
rotation, all the new endpoints in $S_{t+1}$ are at distance two
from a small vertex. It follows by property $(2)$ that $S_{t+1} \cap
SMALL = \emptyset$. Hence $|N(S_{t+1})|\geq 3d|S_{t+1}|$ by property
$(3)$, which implies that after the next rotation we will have
$$\hat{S}_{t+2} \geq \frac{1}{2}(3d |S_{t+1}|- 3(\frac{4}{3}|S_t|
+|S_{t+1}|)) =\frac{3d-7}{2}|S_t|.$$ Hence, a subset
$S_{t+2}\subseteq \hat{S}_{t+2}$ with $|S_{t+2}|=\frac{d}{3}|S_t|$
can be selected.

For the last two rotations our calculations are identical to the ones in
Claim~\ref{St} as those depend on the expansion properties implied by
condition P2.

In conclusion, we created an endpoint set $B(v_1)$ of size at least $n/3$
such that for every $v\in B(v_1)$ there is a $v_1v$-path of maximum length
which can be obtained from $P_0$ by at most
$240+ \frac{4\log n}{\log d}$ rotations with fixed endpoint $v_1$.
{\hfill $\Box$
\medskip\\}

\textbf{Proof of Theorem~\ref{fconnected}}

Let $G=(V,E)$ be $f$-connected where $f(k) = 12 e^{12} + k \log k$.
We prove that $G$ satisfies conditions P1 and P2 with $d=12$ and
apply Theorem~\ref{th} to conclude that $G$ is Hamiltonian for
sufficiently large $n$. Let $A \subseteq V$ be of size at most
$\frac{n}{12 m}$. Either $|A|
> |V \setminus (A \cup N(A))|$ and so in particular $|N(A)| \geq
12|A|$, or the pair $(A \cup N(A),V \setminus A)$ is a separation of
$G$ with $|A| \leq |V \setminus (A \cup N(A))|$ and so by our
assumption $|N(A)| \geq f(|A|) \geq 12 e^{12} + |A| \log |A| \geq
12|A|$. It follows that $G$ satisfies property $P1$ with $d=12$. Let
$A,B$ be two disjoint subsets of $V$ such that $|B| \geq |A| \geq
\frac{n}{4130 m}$. Assume for the sake of contradiction that there
is no edge in $G$ between $A$ and $B$; hence $(V \setminus B, V
\setminus A)$ is a separation of $G$. By our assumption $|(V
\setminus A) \cap (V \setminus B)| = |V \setminus (A \cup B)| \geq
f(|A|) \geq |A| \log |A| > n$. This is clearly a contradiction and
so $G$ satisfies property $P2$ with $d=12$.
{\hfill $\Box$
\medskip\\}

\section*{Acknowledgements} We would like to thank Deryk
Osthus for suggesting to use our criterion to address the conjecture
of~\cite{BBDK}.


\begin{thebibliography}{99}

\bibitem{AKS} M. Ajtai, J. Koml\'{o}s and E. Szemer\'{e}di,
The first occurrence of Hamilton cycles in random graphs,
{\em Annals of Discrete Mathematics} {\bf 27} (1985), 173--178.


\bibitem{Beck}
J. Beck, Random graphs and positional games on the complete graph,
{\em Annals of Discrete Math.} {\bf 28} (1985) 7-13.

\bibitem{Beckbook}
J. Beck, {\bf Tic-Tac-Toe Theory}, {\em manuscript}.

\bibitem{BFF1}
B. Bollob\'as, T. I. Fenner and A. M. Frieze, An algorithm for
finding Hamilton paths and cycles in random graphs, {\em
Combinatorica} {\bf 7} (1987), 327--341.

\bibitem{BBDK}
S. Brandt, H. Broersma, R. Diestel and M. Kriesell, Global
connectivity and expansion: long cycles and factors in $f$-connected
graphs, {\em Combinatorica} {\bf 26} (2006), 17--36.



\bibitem{CE}
V. Chv\'atal and P. Erd\H{o}s, A note on Hamiltonian circuits,
\emph{Discrete Math.} {\bf 2} (1972), 111--113.



\bibitem{Diestel}
R. Diestel, {\bf Graph Theory}, Springer New-York, $2^{nd}$ ed. 1999.



\bibitem{Dirac}
G. Dirac, Some theorems on abstract graphs, \emph{Proc. London Math.
Society} \textbf{2} (1952), 69--81.



\bibitem{ER61} P. Erd\H{o}s and A. R\'enyi, On the evolution of random graphs,
{\em Bull. Inst. Statist. Tokyo} {\bf 38} (1961), 343--347.



\bibitem{FK}
A. Frieze and M. Krivelevich, Hamilton cycles in random
subgraphs of pseudo-random graphs, {\em Discrete Mathematics} {\bf 256} (2002),
137--150.



\bibitem{Gould}
R. J. Gould, Advances on the Hamiltonian problem - a survey,
\emph{Graphs and Combinatorics} {\em 19} (2003), 7--52.



\bibitem{HKS}
D. Hefetz, M. Krivelevich and T. Szab\'o, Avoider-Enforcer games,
{\em Journal of Combinatorial Theory (Series A)}, to appear.



\bibitem{Korsh76} A. D. Korshunov, Solution of a problem of Erd\H{o}s and
R\'enyi on Hamilton cycles in non-oriented graphs, {\em Soviet Math.
Dokl.} {\bf 17} (1976), 760--764.



\bibitem{KS}
J. Koml\'os and E. Szemer\'edi, Limit distributions for the
existence of Hamilton circuits in a random graph, {\em Discrete
Mathematics} {\bf 43} (1983), 55--63.



\bibitem{KrSu}
M. Krivelevich and B. Sudakov, Sparse pseudo-random graphs are
Hamiltonian, {\em Journal of Graph Theory} {\bf 42} (2003), 17--33.



\bibitem{Posa}
L. P\'osa, Hamiltonian circuits in random graphs, {\em Discrete
Math.} {\bf 14} (1976), 359--364.



\bibitem{SV} B. Sudakov and J. Verstra\"ete, Cycle lengths in
sparse graphs, {\em submitted}.



\bibitem{V} J. Verstra\"ete, On arithmetic progressions of
cycle lengths in graphs, {\em Combinatorics, Probability and Computing}
{\bf 9} (2000), 369--373.



\end{thebibliography}
\end{document}